\documentclass{amsart}
\usepackage{amsmath, amstext, amsbsy}

\newtheorem{theorem}{Theorem}[section]

\newtheorem{proposition}{Proposition}[section]   
\newtheorem{prop}{Proposition}[section]

\theoremstyle{definition}

\theoremstyle{remark}

\numberwithin{equation}{section}

\newcommand{\U}{U_q(\widehat{\mathfrak g})}
\renewcommand{\a}{\alpha}
\newcommand{\g}{\gamma}

\newcommand{\la}{\lambda}
\newcommand{\be}{\beta}
\newcommand{\ka}{\kappa}
\newcommand{\ep}{\epsilon}
\newcommand{\vep}{\varepsilon}

\begin{document}


\title[$q$-Vertex operators for quantum affine algebras]
    {$q$-Vertex operators for quantum affine algebras
    }
\author{Naihuan Jing}
\address{Department of Mathematics\\
   North Carolina State University\\
   Raleigh, NC 27695-8205\\
   U.S.A.}
\email{jing@math.ncsu.edu}
\thanks{NJ: Research supported in part by NSA grant MDA904-97-1-0062.}
\author{Kailash C. Misra}
\email{misra@math.ncsu.edu}
\thanks{KCM: Research supported in part by NSA grant
MDA904-96-1-0013.}



\begin{abstract}
$q$-vertex operators for quantum affine algebras have played important role
in the theory of solvable lattice models and the quantum
Knizhnik-Zamolodchikov equation. Explicit constructions of these
vertex operators for most level one modules are known for classical types
except for type $C_n^{(1)}$, where the level $-1/2$ have been
constructed.
In this paper we survey these results for the quantum affine algebras of
types $A_n^{(1)}, B_n^{(1)}, C_n^{(1)}$ and $D_n^{(1)}$.
\end{abstract}

\maketitle

\section{Introduction} \label{S:intro}

Algebraic conformal field theory can be simply described by the
beautiful properties of the Virasoro algebra, vertex operators of affine
Lie algebras, and the Knizhnik-Zamolodchikov equation. The vertex
operators or intertwining operators between certain representations of
affine Lie algebras stand at the core of the foundation since they also
provide solutions to the KZ equation.  Most advances in this part of Lie
theory around 1980s circled around this interesting subject. Vertex
operators also helped the introduction of the vertex algebra by Borcherds
\cite{kn:Bo} (cf. \cite{kn:FLM}).
  In 1985 Drinfeld \cite{kn:Dr1} and Jimbo \cite{kn:Jb}
introduced the notion of quantum group, a
q-deformation of enveloping algebras of Kac-Moody algebras. This has
inspired torrential activity to quantize existed phenomena in Lie
theory. Lusztig first q-deformed the integrable modules of Kac-Moody
algebras \cite{kn:L}. Then the q-deformed level one modules of quantum affine
algebras were constructed by vertex representations
\cite{kn:FJ, kn:J1}. Quantum KZ equations were defined by
Frenkel-Reshetikhin using the q-vertex operators \cite{kn:FR}, and their
matrix coefficients provide the solutions for the q-KZ equations.
 The Kyoto school used the q-vertex operators in solvable
lattice models (see \cite{kn:JM}).
They proposed that the half of the space of states H, on which Baxter's corner
transfer matrix is acting, can be identified with highest weight
representations of quantum affine algebras. Then the correlation
functions and the form factors can be computed in terms of the q-deformed
vertex operators. Moveover their integral formulas are immediately given
by explicit bosonization of the q-vertex operators.
  Because of its applications in
statistical mechanics models and q-conformal field theory,
the program for realizing various q-vertex operators attracted the
attention of several researchers. It started with the
explicit bosonization of level one case for $U_q(A_n^{(1)})$ first for $n=1$
\cite{kn:JMMN},
then for general $n$ by Koyama \cite{kn:Ko1}.
  Kang, Koyama, and the first author \cite{kn:JKK}
generalized the bosonization to type
$D_n^{(1)}$.
More recently we gave explicit bosonizations for other types of quantum
affine algebras \cite{kn:JM1, kn:JM2, kn:JK}. Some lower rank cases
of other modules are also considered in \cite{kn:I, kn:Ko2}.

  In this work we will review the techniques involved in the realization
of q-vertex
operators and provide a unified description for the untwisted
quantum affine
algebras of
classical types. To minimize the length of the paper
we do not include twisted quantum affine algebras which have been
 considered in \cite{kn:JM2}.
However the basic ingredients are similar
in all cases.

Currently two methods of deriving q-vertex operators are used: coproduct
formula (cf. Theorem \ref{T:2.3} ) and direct computation using commutation
relations.
We take the approach of coproduct to present the results, however this
method is not sufficient to determine the q-vertex operator in all the
case (cf. type $C_n^{(1)}$ in \cite{kn:JK}). We refer the reader to
\cite{kn:JK} for more detailed
information about the direct computation approach.
 Due to lack of space we do not include the intertwining operators for the
Drinfeld comultiplication. Their bosonization agrees with our vertex
operators in the special computable component. The rest of the components
in \cite{kn:DI} are then derived using normal ordering and commutation
relations.

 This paper is arranged as follows.
In section two we review the basic notations.
Drinfeld realization, the coproduct
formula and the level zero action of the Drinfeld generators
 are computed in
section three. In section four we discuss the constructions of
level one modules for types $A_n^{(1)}$, $B_n^{(1)}$, $D_n^{(1)}$
and level $-1/2$
modules for type $C_n^{(1)}$. These constructions are
used to construct the q-vertex operators in the last section.


\section{Notation and Preliminaries}

Let
 $A=(A_{ij}), i,j\in I=\{0, 1, \cdots, n\}$ be the  generalized Cartan matrix
of type $A_n^{(1)},
B_n^{(1)}, C_n^{(1)}$, or $D_n^{(1)}$.

Let ${\widehat{\mathfrak g}}={\widehat{\mathfrak g}}(A)$ be the associated
Kac-Moody algebra , and
$h_i$ $(i\in I)$ and $d$ be the generators \cite{kn:K} of the Cartan
subalgebra
 $\widehat {\mathbf h}$.
We define the affine root lattice of $\widehat{\mathfrak g}$ to be
$$\hat{Q}={\mathbb Z}\alpha_0\oplus\cdots\oplus{\mathbb Z}\alpha_n,
$$
where $\alpha_i\in \widehat{\bf h}^*$ such that $<\alpha_i, h_j>=A_{ji}$ and
$<\alpha_i, d>=\delta_{i0}$.
The affine weight lattice $\hat{P}$ is defined to be
$$
\hat{P}={\mathbb Z} \Lambda_0\oplus\cdots\oplus{\mathbb Z}
\Lambda_n\oplus{\mathbb Z}
\delta ,
$$
where
$
\Lambda_i(h_j)=\delta_{ij}, \Lambda_i(d)=0$, and
$\delta(h_j)=0, \delta(d)=1 $ for $j\in I$.
The dual affine weight lattice is then defined as
$$
\hat{P}^{\vee}={\mathbb Z} h_0 \oplus {\mathbb Z} h_1
\oplus \cdots \oplus {\mathbb Z} h_n \oplus {\mathbb Z} d .
$$

We will denote the corresponding weight (resp. root) lattice of
the finite dimensional Lie algebra ${\mathfrak g}$ by $P$ (resp. $Q$).

The nondegenerate symmetric bilinear form $(\ |\ )$ on $\widehat{\bf h}^*$
satisfies
\begin{equation}
(\alpha_i|\alpha_j)=d_iA_{ij}, \ \ (\delta|\alpha_i)=(\delta|\delta)=0
\ \ \mbox{for all} \ i,j\in I \label{E:2.2},
\end{equation}
where $d_i=1$ except for $(d_0, \cdots, d_n)$=$(1, \cdots, 1, 1/2)$ in the case
of $B_n^{(1)}$, and $(d_0, \cdots, d_n)$=$(1, 1/2, \cdots, 1/2, 1)$ in the case
of $C_n^{(1)}$.

Let $q_i=q^{d_i}=q^{\frac 12(\a_i|\a_i)}, i\in I$.
The quantum affine algebra $U_q(\widehat{\mathfrak g})$ is
 the associative algebra with 1 over ${\bf C}(q^{1/2})$
generated by the elements $e_i$, $f_i$ $(i\in I)$ and $q^h$ $(h\in
\hat{P}^{\vee})$
with the following relations :
\begin{eqnarray}
\ q^h q^{h'}&=&q^{h+h'},\quad q^0=1,\ \ \mbox{for} \ h,h'\in \hat{P}^{\vee},
 \nonumber\\
q^h e_{i} q^{-h}&=&q^{\alpha_i(h)} e_{i}, \ \
q^h f_{i} q^{-h}=q^{-\alpha_i(h)} f_{i} \ \ \mbox{for} \ h\in \hat{P}^{\vee}
(i\in I), \nonumber\\
e_{i}f_{j}-f_{j}e_{i} &=&\delta_{ij}
\displaystyle\frac {t_i-t_i^{-1}} {q_i-q_i^{-1}}, \ \mbox{where} \
t_i=q_i^{h_i}=q^{\frac 12(\a_i|\a_i)h_i} \  \mbox{and} \ i,j \in I,
\nonumber
\\
\sum_{m+n=1-a_{ij}}&& (-1)^m e_{i}^{(m)} e_{j} e_{i}^{(n)}=0, \quad
\mbox{and}\nonumber\\
\sum_{m+n=1-a_{ij}}&& (-1)^m f_{i}^{(m)} f_{j} f_{i}^{(n)}=0
 \ \ \mbox{for} \  i\neq j,  \nonumber
\end{eqnarray}
where $e_{i}^{(k)}=e_{i}^k/[k]_i !$, $f_{i}^{(k)}=f_{i}^k/[k]_i !$,
$[m]_i!=\prod_{k=1}^m [k]_i$, and
$[k]_i=\displaystyle\frac {q_i^k-q_i^{-k}} {q_i-q_i^{-1}}.$
For simplicity we will denote $[k]_i=[k]$ for $i=1, \cdots, n-1$.
The derived subalgebra generated by $e_i$, $f_i$, $t_i$ $(i \in I)$ is
denoted by $U'_q(\widehat{\mathfrak g})$.

The algebra $U_q(\widehat{\mathfrak g})$ has a
Hopf algebra structure with comultiplication
$\Delta$, counit $u^*$, and antipode $S$ defined by
\begin{eqnarray}  \nonumber
\Delta(q^h)&=&q^h \otimes q^h \ \ \mbox{for} \ h\in \hat{P}^{\vee},\\
\Delta(e_{i})&=&e_{i}\otimes 1 + t_i \otimes e_{i},
\nonumber\\
\Delta(f_{i})&=&f_{i}\otimes t_i^{-1} + 1 \otimes f_{i} \ \ \mbox{for}
\ i\in I,  \nonumber \\
u^* (q^h)&=&1 \ \ \mbox{for} \ h\in \hat{P}^{\vee}
,\label{E:2.6}\\
u^*(e_{i})&=&u^*(f_{i})=0
\ \ \mbox{for} \ i\in I, \nonumber\\
S(q^h)&=&q^{-h} \ \ \mbox{for} \ h\in \hat{P}^{\vee}, \nonumber\\
S(e_{i})&=&-t_i^{-1}e_{i}, \ \ S (f_{i})=-f_{i}t_i \ \
\mbox{for} \ i\in I.   \nonumber
\end{eqnarray}

Let $V, W$ be two $\U$-modules. The tensor product $V\otimes W$ is defined
as the $\U$-module via the coproduct $\Delta$.
The (restricted) dual $\U$-module $V^*$ is defined by
$$(x \cdot v^*) (u)=v^*(S(x)\cdot u) $$
for $x\in U_q({\widehat{\mathfrak g}})$, $u\in V$, and $v^*\in V^*$.

\section{Drinfeld realization and level zero modules} \label{S:level0}

We now recall Drinfeld's realization of the quantum affine
algebra $U_q(\widehat{\mathfrak g})$ (and of $U_q'(\widehat{\mathfrak g})$)(cf.
\cite{kn:Dr2,kn:B, kn:J3}).
Let ${\bf U}$ be the associative algebra with 1 over ${\bf C}(q^{1/2})$
generated by the elements $x_i^{\pm}(k)$, $a_i(l)$, $K_i^{\pm 1}$,
$\gamma^{\pm 1/2}$, $q^{\pm d}$ $(i=1,2,\cdots,n, k\in {\mathbb Z},
l\in {\mathbb Z} \setminus \{0\})$ with the following defining relations :
\begin{eqnarray}
 [\gamma^{\pm 1/2}, u]&=&0 \ \ \mbox{for all} \ u\in \textstyle {\bf U},
\nonumber\\
\mbox{} [a_i(k), a_j(l)]&=&\delta_{k+l,0}
\displaystyle\frac {[A_{ij}k]_i}{k}
\displaystyle\frac {\gamma^k-\gamma^{-k}}{q_j-q^{-1}_j},
\nonumber\\
\mbox{} [a_i(k), K_j^{\pm 1}]&=&[q^{\pm d}, K_j^{\pm 1}]=0,
\nonumber\\
 q^d x_i^{\pm}(k) q^{-d}&=&q^k x_i^{\pm } (k), \ \
q^d a_i(l) q^{-d}=q^l a_i(l),
\nonumber\\
 K_i x_j^{\pm}(k) K_i^{-1}&=&q^{\pm (\alpha_i|\alpha_j)} x_j ^{\pm}(k),
\nonumber\\
\mbox{} [a_i(k), x_j^{\pm} (l)]&=&\pm \displaystyle\frac {[A_{ij}k]_i}{k}
\gamma^{\mp |k|/2} x_j^{\pm}(k+l),
\nonumber\\
 x_i^{\pm}(k+1)x_j^{\pm}(l)&-&q^{\pm (\alpha_i|\alpha_j)} x_j^{\pm}(l)
x_i^{\pm}(k+1)
\nonumber\\
&=& q^{\pm (\alpha_i|\alpha_j)} x_i^{\pm}(k) x_j^{\pm}(l+1)
-x_j^{\pm}(l+1) x_i^{\pm}(k),
\nonumber\\
\mbox{} [x_i^{+}(k), x_j^{-}(l)]&=&\displaystyle\frac
{\delta_{ij}}{q_i-q_i^{-1}}
\left( \gamma^{\frac {k-l}{2}} \psi_{i} (k+l)-\gamma^{\frac{l-k}{2}}
\varphi_{i} (k+l) \right),
\label{E:2.11}
\end{eqnarray}
where $\psi_{i}(m)$ and $\varphi_{i}(-m)$ $(m\in {\mathbb Z}_{\ge 0})$
are defined by
\begin{eqnarray}
\ && \sum_{m=0}^{\infty} \psi_{i}(m) z^{-m}
=K_i \textstyle {exp} \left( (q_i-q_i^{-1}) \sum_{k=1}^{\infty} a_i(k) z^{-k}
\right),
\nonumber\\
\ && \sum_{m=0}^{\infty} \varphi_{i}(-m) z^{m}
=K_i^{-1} \textstyle {exp} \left(- (q_i-q_i^{-1}) \sum_{k=1}^{\infty} a_i(-k)
z^{k}
\right),
\nonumber
\end{eqnarray}
and the Serre relations are:
\begin{align}
&\mbox{Sym}_{k_1,\cdots,k_m}\sum_{r=0}^{m=1-A_{ij}}(-1)^r
\left[\begin{array}{c} m\\r\end{array}\right]_ix^{\pm}_{i}(k_1)\cdots
x^{\pm}_i(k_r)x^{\pm}_j(l)\\
&\qquad\qquad \cdot x^{\pm}_i(k_{r+1})\cdots x^{\pm}_i(k_m)=0
\hskip 1in\mbox{if} \ i\neq j  .\nonumber
\end{align}
We denote by $\textstyle {\bf U}'$ the subalgebra of {\bf U} generated by
the elements $x_i^{\pm}(k)$, $a_i(l)$, $K_i^{\pm 1}$,
$\gamma^{\pm 1/2}$ $(i=1,2,\cdots,n, k\in {\mathbb Z},
l\in {\mathbb Z} \setminus \{0\})$.

Based on the general result of \cite{kn:Dr2} one can compute directly the
isomorphism as given in \cite{kn:J3}. The $\ep$-sequences correspond
to reduced expressions of the longest element in the Weyl group. For details
see \cite{kn:J3}.
\begin{prop} \cite{kn:J3} \label{P:2.2}
Let $i_1, i_2, \cdots, i_{h-1}$ be a sequence
of indices given in the Table.
 Then the ${\bf C}(q^{1/2})$-algebra isomorphism
$\Psi: U_q(\widehat{\mathfrak g}) \to \textstyle {\bf U}$ is given by
\begin{eqnarray}
e_i &\mapsto& x_i^{+}(0), \ \ f_i \mapsto x_i^{-}(0), \ \
t_i \mapsto K_i \ \ \mbox{for} \ i=1,\cdots, n,
\nonumber\\
e_0 &\mapsto &[x_{i_{h-1}}^{-}(0), \cdots ,x_{i_2}^{-}(0),
x_{i_1}^{-}(1)]_{q^{\ep_1}\cdots q^{\ep_{h-2}}} \g K_{\theta}^{-1},,
\label{E:2.13}\\
f_0 &\mapsto & a(-q)^{-\ep} \g^{-1}K_{\theta}
[x_{i_{h-1}}^{+}(0), \cdots,
x_{i_2}^{+}(0), x_{i_1}^{+}(-1)]_{q^{\ep_1}\cdots q^{\ep_{h-2}}},
\nonumber\\
 t_0 &\mapsto& \gamma K_{\theta}^{-1}, \ \ q^d \mapsto q^d,
\nonumber
\end{eqnarray}
where $K_{\theta}=K_{i_1}\cdots K_{i_{h-1}}$,
$h$ is the Coxeter number, $\ep=\sum_{i=1}^{h-2}\ep_i$.
 The constant $a$ is $1$
for simply laced types $A_n^{(1)}, D_n^{(1)}$, $a=[2]_1$
for $C_n^{(1)}$, and $a=[2]^{1-\delta_{1,i_1}}$ for $B_n^{(1)}$.

The restriction of $\Psi$ to $U_q'(\widehat{\mathfrak g})$ defines an
isomorphism of
$U_q'(\widehat{\mathfrak g})$ and ${\bf U}'$.
\end{prop}

\centerline{Table: $\ep$-Sequences for simple Lie algebras}

\centerline{\vbox{\tabskip=0pt\offinterlineskip
\def\tablerule{\noalign{\hrule}}
\halign to 360pt{\strut#& \vrule#\tabskip=0.1em plus0.1em&
      \hfil#& \vrule#& \hfil#\hfil& \vrule#&
      \hfil#& \vrule#\tabskip=0pt\cr\tablerule
&&\omit\hidewidth $\mathfrak g$\hidewidth&&
  \omit\hidewidth $\ep$-Sequence. $\ep=\sum\ep_i$\hidewidth&&
  \omit\hidewidth $\ep$\hidewidth&\cr\tablerule
&& $A_n$
  && $\a_1 \stackrel{-1}{\rightarrow}\cdots \stackrel{-1}{\rightarrow}
\a_n$
  && -n+1 &\cr\tablerule
&& $B_n$
  && $\a_1 \stackrel{-1}{\rightarrow}\cdots \stackrel{-1}{\rightarrow}
  \a_{n-1} \overset{0}{\rightarrow}\a_n \stackrel{-1}{\rightarrow}
  \cdots \stackrel{-1}{\rightarrow}\a_2$
 && -2n+4 &\cr\tablerule
&& $C_n$
  && $\a_1 \overset{-1/2}{\longrightarrow}\cdots \stackrel{-1}{\rightarrow}
  \a_{n} \overset{-1/2}{\longrightarrow}\cdots \overset{-1/2}{\longrightarrow}
  \a_2 \overset{0}{\rightarrow}\a_1$
 && -n+1 &\cr\tablerule
&& $D_n$
  && $\a_1 \stackrel{-1}{\rightarrow}\cdots \stackrel{-1}{\rightarrow}
  \a_{n} \stackrel{-1}{\rightarrow}\a_{n-2} \stackrel{-1}{\rightarrow}
  \cdots \stackrel{-1}{\rightarrow}\a_2$
 && -2n+4 &\cr\tablerule
}}}

Through the isomorphism of Proposition \ref{P:2.2} the coproduct will be
carried over
to the Drinfeld realization. The following result is true for all type one
quantum affine algebras.
\begin{theorem}\label{T:2.3}
Let $k\in  {\mathbb Z}_{\ge 0}$, $l\in  {\mathbb N}$,
and let $N_{+}^s$ (resp. $N_{-}^s$) be the subalgebra of
{\bf U} generated by the elements $x_{i_1}^{+}(-m_1)\cdots x_{i_s}^{+}(-m_s)$
(resp. $x_{i_1}^{-}(m_1)\cdots x_{i_s}^{-}(m_s)$) with $m_i \in \textstyle
{\mathbb Z}_{\ge 0}$. Then the comultiplication $\Delta$ of the algebra {\bf U}
has the following form:
\begin{eqnarray*}
\Delta(x_{i}^{+}(k))&=&x_i^{+}(k) \otimes \gamma^k
+\gamma^{2k} K_i \otimes x_i^{+}(k)\\
&&\ \ \ \ \ \ \ \ + \sum_{j=0}^{k-1} \gamma^{\frac {k-j}{2}} \psi(k-j)
\otimes \gamma^{k-j} x_i^{+}(j) \ \ (\textstyle {mod} \ N_{-}\otimes N_{+}^2),
\\
 \Delta(x_{i}^{+}(-l))&=&x_i^{+}(-l) \otimes \gamma^{-l}
+ K_i^{-1} \otimes x_i^{+}(-l)\\
&&\ \ \ \ \ \ \ \ + \sum_{j=1}^{l-1} \gamma^{\frac {l-j}{2}} \varphi(-l+j)
\otimes \gamma^{-l+j} x_i^{+}(-j) \ \ (\textstyle {mod} \ N_{-}\otimes
N_{+}^2), \\
 \Delta(x_{i}^{-}(l))&=&x_i^{-}(l) \otimes K_i
+\gamma^{l} \otimes x_i^{-}(l)\\
&&\ \ \ \ \ \ \ \ + \sum_{j=1}^{l-1} \gamma^{l-j} x_i^{-}(j)
\otimes \gamma^{\frac {j-l}{2}} \psi_i(l-j)
 \ \ (\textstyle {mod} \ N_{-}^2 \otimes N_{+}), \\
 \Delta(x_{i}^{-}(-k))&=&x_i^{-}(k) \otimes \gamma^{-2k}K_i^{-1}
+\gamma^{-k} \otimes x_i^{-}(-k)\\
&&\ \ \ \ \ \ \ \ + \sum_{j=0}^{k-1} \gamma^{j-k} x_i^{-}(-j)
\otimes \gamma^{-\frac {k+3j}{2}} \varphi_i(j-k)
 \ \ (\textstyle {mod} \ N_{-}^2 \otimes N_{+}), \\
 \Delta(a_i(l))&=&a_i(l) \otimes \gamma^{\frac {l}{2}}
+ \gamma^{\frac {3l}{2}} \otimes a_i(l) \ \ (\textstyle {mod} \ N_{-}\otimes
N_{+}),\\
 \Delta(a_i(-l))&=&a_i(-l) \otimes \gamma^{-\frac {3l}{2}}
+ \gamma^{-\frac {l}{2}} \otimes a_i(-l) \ \ (\textstyle {mod} \ N_{-}\otimes
N_{+}).
\end{eqnarray*}
Moreover the same formulas are true for the derived subalgebra
$\textstyle {\bf U}'$.
\end{theorem}

Let $V$ be the finite dimensional representation described by the following
representation graphs. Note that the graphs are crystal graphs, and
they are also perfect crystals \cite{kn:KMN}
except in the case of $C_n^{(1)}$.
Let $v_i$ be the basis elements ($i\in I$), and $I$ is
the index set. Let $E_{ij}$ be the matrix units. We can read the actions from
the graphs.

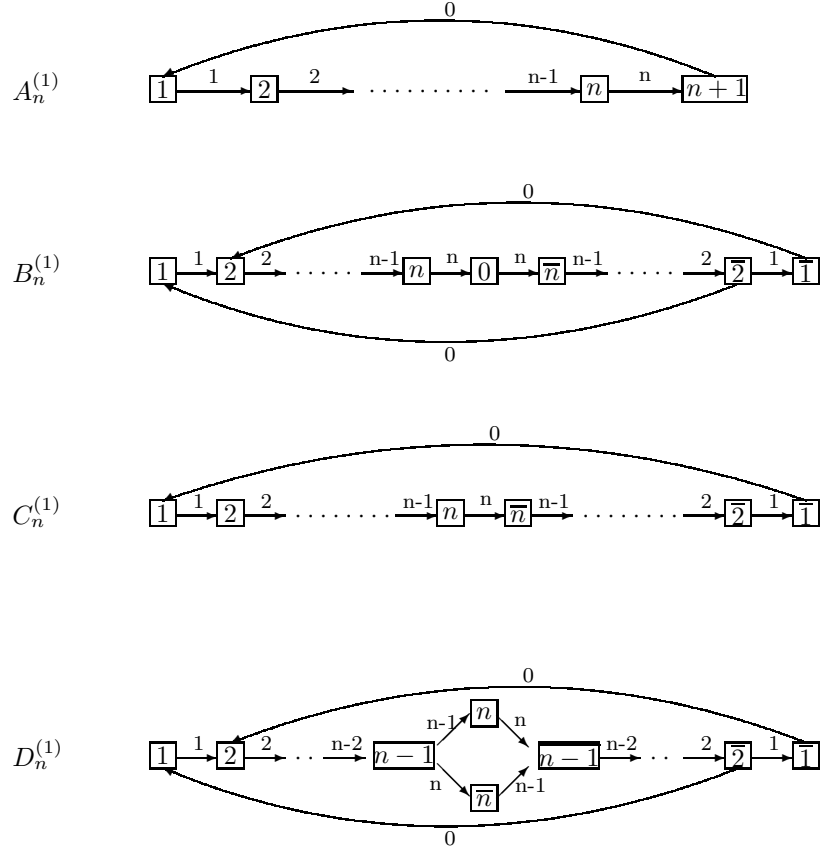
\begin{figure}[ht]
\setlength{\unitlength}{0.8pt}
\hspace{1cm}
\begin{picture}(300,60)(-50,0)
\put(-60,16){$A_n^{(1)}$}
\put(5, 15){\framebox(12, 12){$1$}}
\put(17,20){\vector(1,0){36}}
\put(53, 15){\framebox(12, 12){$2$}}
\put(65,20){\vector(1,0){36}}
\multiput(110,20)(6,0){10}{\circle*{1}}
\put(173,20){\vector(1,0){36}}
\put(209, 15){\framebox(12, 12){$n$}}
\put(221, 20){\vector(1,0){36}}
\put(257, 15){\framebox(30, 12){$\, n+1$}}
\qbezier(272,27)(142, 80)(11, 27)    
\put(13,28){\vector(-2,-1){2}}
\put(30, 22){\makebox(10,10){\footnotesize 1}}
\put(78, 22){\makebox(10,10){\footnotesize 2}}
\put(186, 22){\makebox(10,10){\footnotesize n-1}}
\put(234,22){\makebox(10,10){\footnotesize n}}
\put(142, 54){\makebox(10,10){\footnotesize 0}}
\end{picture}
\vskip0.7cm
\hspace{1cm}
\begin{picture}(300,60)(-50,0)
\put(-60,16){$B_n^{(1)}$}
\put(5, 15){\framebox(12, 12){$1$}}
\put(17,20){\vector(1,0){20}}
\put(37, 15){\framebox(12, 12){$2$}}
\put(49,20){\vector(1,0){20}}
\multiput(75,20)(6,0){5}{\circle*{1}}
\put(105,20){\vector(1,0){20}}
\put(125, 15){\framebox(12, 12){$n$}}
\put(137, 20){\vector(1,0){20}}
\put(157, 15){\framebox(12, 12){$0$}}
\put(169, 20){\vector(1,0){20}}
\put(189, 15){\framebox(12,12){$\overline n$}}
\put(201,20){\vector(1,0){20}}
\multiput(227,20)(6,0){5}{\circle*{1}}
\put(257, 20){\vector(1,0){20}}
\put(277, 15){\framebox(12,12){$\overline 2$}}
\qbezier(283,15)(142, -40)(11, 15)    
\put(13, 14){\vector(-2, 1){2}}
\put(289,20){\vector(1,0){20}}
\put(309, 15){\framebox(12,12){$\overline 1$}}
\qbezier(315,27)(179, 80)(43, 27)    
\put(45, 28){\vector(-2, -1){2}}
\put(23, 22){\makebox(10,10){\footnotesize 1}}    
\put(55, 22){\makebox(10,10){\footnotesize 2}}
\put(111, 22){\makebox(10,10){\footnotesize n-1}}
\put(143,22){\makebox(10,10){\footnotesize n}}
\put(175,22){\makebox(10,10){\footnotesize n}}
\put(207, 22){\makebox(10,10){\footnotesize n-1}}
\put(263, 22){\makebox(10,10){\footnotesize 2}}
\put(295, 22){\makebox(10,10){\footnotesize 1}}
\put(142, -23){\makebox(10,10){\footnotesize 0}}
\put(179, 54){\makebox(10,10){\footnotesize 0}}
\end{picture}
\vskip 1.5cm
\hspace{1cm}
\begin{picture}(300,60)(-50,0)
\put(-60,16){$C_n^{(1)}$}
\put(5, 15){\framebox(12, 12){$1$}}
\put(17,20){\vector(1,0){20}}
\put(37, 15){\framebox(12, 12){$2$}}
\put(49,20){\vector(1,0){20}}
\multiput(75,20)(6,0){8}{\circle*{1}}
\put(121,20){\vector(1,0){20}}
\put(141, 15){\framebox(12, 12){$n$}}
\put(153, 20){\vector(1,0){20}}
\put(173, 15){\framebox(12,12){$\overline n$}}
\put(185,20){\vector(1,0){20}}
\multiput(211,20)(6,0){8}{\circle*{1}}
\put(257, 20){\vector(1,0){20}}
\put(277, 15){\framebox(12,12){$\overline 2$}}
\put(289,20){\vector(1,0){20}}
\put(309, 15){\framebox(12,12){$\overline 1$}}
\qbezier(315,27)(163, 80)(11, 27)    
\put(13, 28){\vector(-2, -1){2}}
\put(23, 22){\makebox(10,10){\footnotesize 1}}    
\put(55, 22){\makebox(10,10){\footnotesize 2}}
\put(127, 22){\makebox(10,10){\footnotesize n-1}}
\put(159,22){\makebox(10,10){\footnotesize n}}
\put(191, 22){\makebox(10,10){\footnotesize n-1}}
\put(263, 22){\makebox(10,10){\footnotesize 2}}
\put(295, 22){\makebox(10,10){\footnotesize 1}}
\put(163, 54){\makebox(10,10){\footnotesize 0}}
\end{picture}
\vskip1.5cm
\hspace{1cm}
\begin{picture}(300,60)(-50,0)
\put(-60,16){$D_n^{(1)}$}
\put(5, 15){\framebox(12, 12){$1$}}
\put(17,20){\vector(1,0){20}}
\put(37, 15){\framebox(12, 12){$2$}}
\put(49,20){\vector(1,0){20}}
\multiput(75,20)(6,0){2}{\circle*{1}}
\put(87,20){\vector(1,0){20}}
\put(111, 15){\framebox(28, 12){$n-1$}}
\put(141, 26){\vector(1,1){15}}
\put(157, 35){\framebox(12, 12){$n$}}
\put(169, 40){\vector(1,-1){15}}
\put(141, 17){\vector(1,-1){15}}  %
\put(157, -5){\framebox(12, 12){$\overline n$}}
\put(169, 1){\vector(1, 1){15}}   %
\put(189, 15){\framebox(28,12){$\overline {n-1}$}}
\put(217,20){\vector(1,0){20}}
\multiput(243,20)(6,0){2}{\circle*{1}}
\put(257, 20){\vector(1,0){20}}
\put(277, 15){\framebox(12,12){$\overline 2$}}
\qbezier(283,15)(142, -40)(11, 15)    
\put(13, 14){\vector(-2, 1){2}}
\put(289,20){\vector(1,0){20}}
\put(309, 15){\framebox(12,12){$\overline 1$}}
\qbezier(315,27)(179, 80)(43, 27)    
\put(45, 28){\vector(-2, -1){2}}
\put(23, 22){\makebox(10,10){\footnotesize 1}}    
\put(55, 22){\makebox(10,10){\footnotesize 2}}
\put(93, 22){\makebox(10,10){\footnotesize n-2}}
\put(138,32){\makebox(10,10){\footnotesize n-1}}
\put(175,32){\makebox(10,10){\footnotesize n}}
\put(135,2){\makebox(10,10){\footnotesize n}}
\put(180,2){\makebox(10,10){\footnotesize n-1}}
\put(223, 22){\makebox(10,10){\footnotesize n-2}}
\put(263, 22){\makebox(10,10){\footnotesize 2}}
\put(295, 22){\makebox(10,10){\footnotesize 1}}
\put(142, -23){\makebox(10,10){\footnotesize 0}}
\put(179, 54){\makebox(10,10){\footnotesize 0}}
\end{picture}
\vskip 0.5cm
\caption{Representation graphs}
\end{figure}

The $\U$-module structure on the affinization or evaluation
module of $V$ \cite{kn:KMN} is given as follows.
We equip the affinization $V_{z}=V \otimes {\bf C}(q^{1/2})[z, z^{-1}]$
with the following
actions:
\begin{eqnarray}
 e_i(v\otimes z^m)&=&e_i v \otimes z^{m+\delta_{i,0}}, \ \
 f_i(v\otimes z^m)=f_i v \otimes z^{m-\delta_{i,0}},
 \nonumber\\
 t_i(v\otimes z^m)&=&t_i v \otimes z^{m},  \label{E:2.18}\\
 q^d (v\otimes z^m)&=& q^m v\otimes z^m,   \nonumber
\end{eqnarray}
for $i\in I$, $v\in V$.

The evaluation module $V_z$ is a level zero $\U$-module, i.e., $\gamma$ acts
as identity ($=q^0$). Through the isomorphism $\Psi$ the evaluation module
is also a $\bf U$-module. The action of the Drinfeld generators are given
by the following result.

Define $\tau=\tau(A)$ as follows:
\begin{equation*}
\tau=\left\{\begin{array}{ll}  n+1 & \mbox{$A_n^{(1)}$}\\
2n-1 & \mbox{$B_n^{(1)}$}\\
2n+2  & \mbox{$C_n^{(1)}$}\\
2n-2 & \mbox{$D_n^{(1)}$.}\end{array}\right.
\end{equation*}
We define the content of each vertex in the representation graph by
\begin{equation*}
c(i)=\left\{\begin{array}{ll}  i & \mbox{ if $i < n$}\\
\tau -\overline{i}  & \mbox{if $i \in \{\overline 1, \cdots, \overline n\}$}
\end{array}\right.
\end{equation*}
Here we adopt the convention that $\overline{\overline i}=i$.

\begin{theorem}\cite{kn:Ko1, kn:JKK, kn:JM, kn:JK}
\label{T:2.4} For $j\in \{1, \cdots , n\}$,
the Drinfeld generators act on the
evaluation module $V_z$ as follows.
\begin{equation}\label{E:2.20}
\begin{array}{rcl}
x_j^{+}(k)&=&\sum_{(i(j), e(j))} (q^{c(j)}_jz)^k E_{i(j), e(j)}\\
x_j^{-}(k)&=&\sum_{(i(j), e(j))} (q^{c(j)}_jz)^k E_{e(j), i(j)}\\
a_j(l)&=&\sum_{(i(j), e(j))} \frac{[l]}l (q^{c(j)}_jz)^l
(q^{-l}_jE_{i(j), i(j)}-q^l_jE_{e(j), e(j)})\\
x_n^{+}(k)&=&\left\{\begin{array}{ll}
(q^{n}z)^k[2]_nE_{n,0}+
(q^{n-1}z)^kE_{0, \overline {n}}, & \mbox{ $B_n^{(1)}$}\\
(q^{n-1}z)^kE_{n,\overline n}, & \mbox{ $C_n^{(1)}$}\\
(q^{n-1}z)^k(E_{n-1, \overline n}+E_{n, \overline{n-1}}), & \mbox{ $D_n^{(1)}$}
\end{array}\right.\\
x_n^{-}(k)&=&\left\{\begin{array}{ll}
(q^{n-1}z)^k[2]_nE_{\overline n, 0}+(q^{n}z)^k
E_{0,n} & \mbox{ $B_n^{(1)}$}\\
(q^{n-1}z)^kE_{\overline n, n}, & \mbox{ $C_n^{(1)}$}\\
(q^{n-1}z)^k(E_{\overline n, n-1}+E_{\overline{n-1}, n}), & \mbox{ $D_n^{(1)}$}
\end{array}\right.\\
a_n(l)&=&\left\{\begin{array}{ll}
\frac {[2l]_n}{l}  (q^{n-1}z)^l ( (E_{n,n}
-q^l E_{0, 0}) + (E_{0, 0}
-q^l E_{\overline {n}, \overline {n}}) ) & \mbox{ $B_n^{(1)}$}\\
\frac{[l]}l(q^{n-1}z)^l(q^{-l}E_{nn}-q^lE_{\overline n, \overline n}),
& \mbox{ $C_n^{(1)}$}\\
\frac {[l]}{l}  (q^{n-1}z)^l ( q^{-l}E_{n-1,n-1}
-q^l E_{\overline n, \overline n} &\\
 \qquad\qquad + q^{-l}E_{n, n}
-q^l E_{\overline {n-1}, \overline {n-1}} ) & \mbox{ $D_n^{(1)}$}
\end{array}\right.
\end{array}
\end{equation}
where the summation runs through all possible pairs
$i(j)\stackrel{j}{\longrightarrow}
e(j)$ in the representation graphs connected by $j$.
$j=1, 2, \cdots, n-1$
(but we allow $j=n$ for type $A_n^{(1)}$), $k\in {\mathbb Z}$ and $l\in
{\mathbb Z}
\setminus \{0\}$.
\end{theorem}

\begin{theorem} \cite{kn:Ko1, kn:JKK, kn:JM, kn:JK}\label{T:2.5} The Drinfeld
generators act on the
dual evaluation module $V_z^*=V^*\otimes {\bf C}(q^{1/2})[z, z^{-1}]$
as follows.
\begin{equation}\label{E:2.21}
\begin{array}{rcl}
x_j^{+}(k)&=&(-q^{-1})\sum_{(i(j), e(j))} (q^{-c(j)}_jz)^k E_{e(j), i(j)}^*\\
x_j^{-}(k)&=&(-q)\sum_{(i(j), e(j))} (q^{-c(j)}_jz)^k E_{i(j), e(j)}^*\\
a_j(l)&=&\sum_{(i(j), e(j))} \frac{[l]}l (q^{-c(j)}_jz)^l
(q^{-l}_jE_{e(j), e(j)}^*-q^l_jE_{i(j), i(j)}^*)\\
x_n^{+}(k)&=&\left\{\begin{array}{ll}
(-q^{-1})((q^{-n}z)^k[2]_nE_{0, n}^*+
(q^{-n+1}z)^kqE_{\overline {n}, 0}^*), & \mbox{ $B_n^{(1)}$}\\
(-q^{-1})(q^{-n+1}z)^kE_{\overline n, n}^*, & \mbox{ $C_n^{(1)}$}\\
(-q^{-1})(q^{-n+1}z)^k(E_{ \overline n, n-1}^*+E_{\overline{n-1}, n}^*),
& \mbox{ $D_n^{(1)}$}
\end{array}\right.\\
x_n^{-}(k)&=&\left\{\begin{array}{ll}
(-q)((q^{-n+1}z)^k[2]_nq^{-1}E_{0, \overline n}^*+(q^{-n}z)^k
E_{n, 0}^*) & \mbox{ $B_n^{(1)}$}\\
(-q)(q^{-n+1}z)^kE_{n, \overline n}^*, & \mbox{ $C_n^{(1)}$}\\
(-q)(q^{-n+1}z)^k(E_{n-1, \overline n}^*+E_{n, \overline{n-1}}^*),
& \mbox{ $D_n^{(1)}$}
\end{array}\right.\\
a_n(l)&=&\left\{\begin{array}{ll}
\frac {[2l]_n}{l} (q^{-n+1}z)^l
((q^{-l} E_{0, 0}^*-  E_{n,n}^*)
+ (q^{-l} E_{\overline {n}, \overline {n}}^*-E_{0, 0}^*
) ), & \mbox{$B_n^{(1)}$}\\
\frac{[l]}l(q^{-n+1}z)^l(q^{-l}E_{\overline n, \overline n}^*-q^{l}E_{nn}^*),
& \mbox{ $C_n^{(1)}$}\\
\frac {[l]}{l}  (q^{-n+1}z)^l ( q^{-l} E_{\overline n, \overline
n}^*-q^{l}E_{n-1,n-1}^*
&\\
 \qquad\qquad + q^{-l} E_{\overline {n-1}, \overline {n-1}}^*-q^{l}E_{n, n}^*
 ) & \mbox{ $D_n^{(1)}$}
\end{array}\right.
\end{array}
\end{equation}
where the summation runs through all possible pairs
$i(j)\stackrel{j}{\longrightarrow}
e(j)$ in the representation graphs connected by $j$.
$j=1, 2, \cdots, n-1$
(but we allow $j=n$ for type $A_n^{(1)}$), $k\in {\mathbb Z}$ and $l\in
{\mathbb Z}
\setminus \{0\}$.
\end{theorem}

\section{Vertex representation of $\U$} \label{S:level1}

In this section we recall the level one realizations of $\U$ for
$\widehat{\mathfrak g}=A_n^{(1)}, D_n^{(1)}$ as given in \cite{kn:FJ}
and $\widehat{\mathfrak g}=B_n^{(1)}$ in \cite{kn:JM1} (cf. \cite{kn:Be} for
$V(\Lambda_0)$). We also recall
the level $-1/2$ realization of $U_q(C_n^{(1)})$ as given in \cite{kn:JKM}.
The quantum affine algebra $\U$ has level one
irreducible representations:
$V(\Lambda_i)$, $i=0, 1, \cdots, n$ for $A_n^{(1)}$;
$i=0, 1, n$ for $B_n^{(1)}$ and $i=0, 1, n-1, n$ for $D_n^{(1)}$.

The weight lattice $P$ has two coset representatives $Q$ and $Q+\lambda_n$,
where $\la_i$ denote the fundamental weights in $P$. We will express
any level one fundamental weight $\la_i$ as a linear combination of
the simple roots $\a_i$ and $\la_n$.

Let $\vep$: $Q\times Q\longrightarrow {\mathbb Z}_2=\{\pm 1\}$ be a
cocycle such that
\begin{equation}\label{E:3.0}
\vep(\a, \be)\vep(\be, \a)=(-1)^{(\a|\be)+(\a|\a)(\be|\be)},
\end{equation}
whose existence can be constructed directly on the simple roots.
Then there corresponds a central extension of the group algebra
${\bf C}(Q)$ associated with $\vep$ such that
$$
1\longrightarrow {\mathbb Z}_2\longrightarrow {\bf C}\{Q\}
\longrightarrow  {\bf C}(Q) \longrightarrow 1.
$$
If we still use $e^{\a}$
to denote the generators of ${\bf C}\{Q\}$, then we
have
\begin{equation}\label{E:3.3}
e^{\a}e^{\be}=(-1)^{(\a|\be)+(\a|\a)(\be|\be)}e^{\be}e^{\a}.
\end{equation}
for any $\a, \be\in Q$.

Moreover we define the vector space ${\bf C}\{P\}={\bf C}\{Q\}\oplus
{\bf C}\{Q\}e^{\la_n}$ as a
${\bf C}\{Q\}$-module by formally adjoining the element
$e^{\la_n}$.

For $\a\in P$ define the operator $\partial_{\alpha}$ on {\bf C}\{P\}
by
\begin{equation}\label{E:3.5}
\partial_{\alpha}e^{\be}=(\a|\be)e^{\be}
\end{equation}

Let $U_q(\hat{\bf h})$ be the infinite dimensional Heisenberg algebra
generated by $a_i(k)$ and the central element
$\gamma$ ($k\in {\mathbb Z}\backslash\{0\}, i=1, \cdots,
n$) subject to the relations
\begin{equation}\label{E:3.6}
[a_i(k), a_j(l)]=\delta_{k+l,0}
\displaystyle\frac {[A_{ij}k]_i}{k}
\displaystyle\frac {\gamma^k-\gamma^{-k}}{q_j-q^{-1}_j}.
\end{equation}
The algebra $U_q(\hat{\bf h})$ acts on the space of symmetric algebra
$Sym(\hat{\bf h}^-)$ generated by
$a_i(-k)$ ($k\in {\mathbb N}, i=1, \cdots, n$) with $\gamma =q$
by the action:
\begin{eqnarray*}
&a_i(-k)&\mapsto \mbox{multiplication by}\ a_i(-k),\\
&a_i(k)&\mapsto \sum_{j}\frac{[A_{ij}k]_i}{k}\frac{q^k-q^{-k}}{q_j-q_j^{-1}}
\frac{d}{d\,a_j(-k)}.
\end{eqnarray*}

For the case of $B_n^{(1)}$ we need to consider an extra
fermionic field. For $Z={\mathbb Z}+s$ ($s=1/2 $, NS-case; $s =0$, R-case), let
$C_q^{s}$ be the
$q$-deformed Clifford
algebra generated by $\kappa(k)$, $k\in Z$
satisfying the relations:
\begin{equation}\label{E:3.7}
\{\kappa(k), \ka(l)\}=(q^{k}+q^{-k})\delta_{k, -l},
\end{equation}
where $k, l \in Z$ and $\{a, b\}=ab+ba$.

Let $\Lambda(C^-_q)^s$ be the polynomial algebra generated by $\ka(-k)$,
$k\in Z_{>0}$, and $\Lambda(C^-_q)^s_0$ (resp. $\Lambda(C^-_q)^s_1$) be the
subalgebra
generated by products of even (resp. odd) number of generators $\ka(-k)$'s.
Then
$$
\Lambda(C^-_q)^s=\Lambda(C^-_q)^s_0\oplus \Lambda(C^-_q)^s_1.
$$

The algebra $C_q^{s}$ acts on $\Lambda(C^-_q)^s$ canonically by
($k\in Z_{>0}$)
\begin{eqnarray*}
&\ka(-k)&\mapsto \mbox{multiplication by}\ \ka(-k),\\
&\ka(k)&\mapsto -\frac{d}{d\ \ka({-k})}.
\end{eqnarray*}

We define that
\begin{eqnarray}\nonumber
V(\Lambda_i)&=&Sym\hat{\bf h}^-)\otimes ({\bf C}\{Q\}e^{\la_i},
\qquad\mbox{$A_n^{(1)}, D_n^{(1)}$}\\
V(\Lambda_0)&=&Sym(\hat{\bf h}^-)\otimes ({\bf C}\{Q_0\}
\otimes\Lambda(C^-_q)^{1/2}_0\oplus {\bf C}\{Q_0\}e^{\lambda_1}
\otimes\Lambda(C^-_q)^{1/2}_1), \mbox{$B_n^{(1)}$}\nonumber \\
V(\Lambda_1)&=&Sym(\hat{\bf h}^-)\otimes ({\bf C}\{Q_0\}e^{\lambda_1}
\otimes\Lambda(C^-_q)^{1/2}_0\oplus {\bf C}\{Q_0\}
\otimes\Lambda(C^-_q)^{1/2}_1),
\mbox{$B_n^{(1)}$}\nonumber\\
V(\Lambda_n)&=&Sym(\hat{\bf h}^-)\otimes {\bf C}\{Q\}e^{\la_n}\otimes
\Lambda(C^-_q)^0, \mbox{$B_n^{(1)}$}
\nonumber
\end{eqnarray}
where $i=0, \cdots, n$ for $A_n^{(1)}$ and $i=0, 1, n-1, n$ for
$D_n^{(1)}$. ${\bf C}\{Q_0\}$ is the group algebra of the sublattice of
long roots in the case of $B_n^{(1)}$. Here for convenience we
denote $\la_0=0$.

\begin{prop}\cite{kn:FJ}\cite{kn:JM1}\label{P:3.1}
The space $V(\Lambda_i)$ of level one irreducible modules
of the quantum affine Lie algebra
$U_q(\widehat{\mathfrak g})$ is realized under the following action:
\begin{eqnarray*}
 \gamma &\mapsto & q, \
K_j \mapsto q^{\partial_{\a_j}},  \  a_j(k)\mapsto a_j(k)
\quad (1\leq j\leq n),\\
x_i^{\pm}(z)&\mapsto& X^{\pm}_i(z)
=exp(\pm\sum_{k=1}^{\infty}\frac{a_i(-k)}{[k]}q^{\mp k/2}z^k)
 exp(\mp\sum_{k=1}^{\infty}\frac{a_i(k)}{[k]}q^{\mp k/2}z^{-k})\\
& &\times e^{\pm \a_i}z^{\pm \partial_{\alpha_i}+1},
    \qquad\mbox{$1\leq i\leq n$;
when $\widehat g=B_n^{(1)}$, $1\leq i\leq n-1$} \\
 x_n^{\pm}(z)&\mapsto& X^{\pm}_n(z)
=exp(\pm\sum_{k=1}^{\infty}\frac{a_n(-k)}{[2k]_n}q^{\mp k/2}z^k)
 exp(\mp\sum_{k=1}^{\infty}\frac{a_n(k)}{[2k]_n}q^{\mp k/2}z^{-k})\\
& &\times e^{\pm \a_n}z^{\pm \partial_{\alpha_n}+1/2}(\pm\ka(z)) ,
\qquad\mbox{
for $\widehat g=B_n^{(1)}$}
\end{eqnarray*}
where  $\ka(z)=\sum_{m\in Z}\ka(m)z^{-m}$.
The highest weight vectors are
$
 |\Lambda_i\rangle
=e^{\lambda_i}.
$
\end{prop}

 We now recall the level $-1/2$ realization of $U_q(C_n^{(1)})$ \cite{kn:JKM}.
  There are four level $-1/2$ weights:
$\mu_1=-\frac 12\Lambda_0, \mu_2= -\frac 32\Lambda_0+\Lambda_1-\frac 12\delta$,
$\mu_3= -\frac 12\Lambda_n, \mu_4=-\frac 32\Lambda_n+\Lambda_{n-1}$.

Let $a_i(m)$ be the bosonic operators satisfying (\ref{E:3.6}) with
$\gamma=q^{-1/2}$.
Let $b_i(m)$ be another set of bosonic operators
   satisfying the following defining relations.

   \begin{equation}
   \begin{array}{rcl}
   \ [b_i(m), b_j(l)]&=&m\delta_{ij}\delta_{m+l,0},\\
   \ [a_i(m), b_j(l)]&=&0 .
   \end{array}
   \end{equation}
   We define the Fock space
   ${\mathcal F}_{\a, \beta}$
   for $\a \in P+\frac{\mathbb Z}{2}\la_n$, $\be \in P$
   by the defining relations
   $$ a_i(m) | \a, \be \rangle = 0 \quad (m>0)
      \; , \quad
      b_i(m) | \a, \be \rangle = 0 \quad (m>0) \;, $$
   $$ a_i(0) | \a, \be \rangle = (\a_i|\a) | \a, \be \rangle
      \; , \quad
      b_i(0) | \a, \be \rangle = (2\vep_i|\be) | \a, \be \rangle \; , $$
   where $| \a, \be \rangle$ is the vacuum vector, and $\vep_i$
are the usual orthonormal vectors.

   We set
   $$ \widetilde{\mathcal F}=
	  \bigoplus_{\a\in P+\frac12{\mathbb Z}\la_n,\be \in P}
	  {\mathcal F}_{\a, \be} . $$
   Let $e^{a_i}=e^{\a_i}$ and $e^{b_i}$ ($1\leq i\leq n$) be operators
   on ${\widetilde{\mathcal F}}$ given by:
   $$ e^{\a_i}|\a,\be \rangle= |\a+\vep_i,\be \rangle
      \quad , \quad
      e^{b_i}|\a, \be \rangle= |\a,\be+\vep_i \rangle . $$
   Let $: \quad :$ be the usual bosonic normal ordering defined by
   $$ :a_i(m) a_j(l): = a_i(m) a_j(l) \, (m \leq l) ,
   \; a_j(l) a_i(m) \, (m>l) , $$
   $$ :e^{a} a_i(0):=:a_i(0) e^{a}:=e^{a} a_i(0) \, . $$
   and similar normal products for the $b_i(m)$'s.
   Let $\partial=\partial_{q^{1/2}}$ be the $q$-difference operator:
   $$ \partial_{q^{1/2}}(f(z))
   =\frac{f(q^{1/2}z)-f(q^{-1/2}z)}{(q^{1/2}-q^{-1/2})z} $$
   We introduce the following operators.
\begin{eqnarray*}
&Y_i^{\pm}(z)&=
	 \exp ( \pm \sum^{\infty}_{k=1}
		 \frac{a_i(k)}{[-\frac{1}{2}k]} q^{\pm \frac{k}{4}} z^k)
	 \exp ( \mp \sum^{\infty}_{k=1}
	   \frac{a_i(k)   }{[-\frac{1}{2}k]} q^{\pm \frac{k}{4}} z^{-k})
	 e^{\pm a_i} z^{\mp 2a_i(0)} ,\\
&Z^{\pm}_i(z)&=
	     \exp ( \pm \sum^{\infty}_{k=1} \frac{b_i(-k)}{k} z^k)
	     \exp ( \mp \sum^{\infty}_{k=1} \frac{b_i(k)}{k} z^{-k})
	     e^{\pm b_i} z^{\pm b_i(0)} .
\end{eqnarray*}
   We define the operators $x_i^{\pm}(m)$ $(i=1, \cdots, n, m \in {\mathbb Z})$
   by the following generating functions.
   $X_i^{\pm}(z)=\sum_{m \in {\mathbb Z}} x_i^{\pm}(m) z^{-m-1}$.
\begin{eqnarray*}
X_i^+(z) &=&
	     \partial Z^+_i(z)Z^-_{i+1}(z)Y_i^+(z), \qquad i=1, \cdots n-1\\
X_i^-(z) &=&
	     Z^-_i(z)\partial Z^+_{i+1}(z)Y^-_i(z), \qquad i=1, \cdots n-1\\
X_n^+(z) &=&
	     \left(
	     \frac{1}{q^{\frac{1}{2}}+ q^{-\frac{1}{2}}}
	     :Z^+_n(z) \partial^2 Z^+_n(z):
     -:\partial Z^+_n(q^{\frac{1}{2}}z)
	      \partial Z^+_n(q^{-\frac{1}{2}}z):
	     \right)
	     Y^+_n(z) \\
X_n^-(z) &=&
	     \frac{1}
		  {q^{\frac{1}{2}}+ q^{-\frac{1}{2}}}
	     :Z^-_n(q^{\frac{1}{2}}z)Z^-_n(q^{-\frac{1}{2}}z): Y^-_n(z).
\end{eqnarray*}

\noindent{\it Remark}. Our $a_i(k)$ differs from that in \cite{kn:JKM}, where
we took $a_i(k)/[d_i]$ for $a_i(k)$.

  Furthermore, we know that $\widetilde{\mathcal F}$ contains
  the four irreducible highest weight modules (\cite{kn:JKM}).
  Let ${\mathcal F}_{\a,\be}^1$ be the subspace of ${\mathcal F}_{\a,\be}$
  generated by $a_i(m)$ $(i=1,\cdots,n\in \mathbb Z$, $m\in \mathbb Z$)
Similarly, let ${\mathcal F}_{\a,\be}^{2,j}$ $(j=1, \cdots, n)$ be
  the subspace of ${\mathcal F}_{\a,\be}$ generated by
  $b_j(m)$ $(m \in {\mathbb Z})$.
  We can define the following isomorphism by
  $|\a,\be\rangle \otimes |\a',\be' \rangle
   \rightarrow |\a+\a',\be+\be' \rangle$.
  $$ {\mathcal F}^1_{\a,0} \otimes
     {\mathcal F}^{2,1}_{0,\be_1} \otimes \cdots \otimes {\mathcal
F}^{2,n}_{0,\be_n}
     \longrightarrow
     {\mathcal F}_{\a, \be_1+\cdots+\be_n} . $$
  Let $Q_j^-$ be the operator from ${\mathcal F}^{2,j}_{\a, \be}$
  to ${\mathcal F}^{2,j}_{\a, \be-\vep_j}$ defined by
  \[ Q_i^-=\frac{1}{2 \pi \sqrt{-1}} \oint Z^-_i(z) dz .\]
  We set subspaces ${\mathcal F}_i$ $(i=1,2,3,4)$
  of $\widetilde{\mathcal F}$ as follows.

\[\begin{array}{llll}
  \displaystyle{
 {\mathcal F}_1 = \bigoplus_{\alpha \in Q}
		   {\mathcal F}'_{\alpha,\alpha}},
    &&\displaystyle{{\mathcal F}_2 = \bigoplus_{\alpha \in Q}
		   {\mathcal F}'_{\alpha+\vep_1,\alpha+\vep_1}}\\
    \displaystyle{{\mathcal F}_3 = \bigoplus_{\alpha \in Q}
		   {\mathcal F}'_{\alpha-\frac{1}{2}\lambda_n, \alpha}},
   &&\displaystyle{{\mathcal F}_4 = \bigoplus_{\alpha \in Q +\vep_n}
		   {\mathcal F}'_{\alpha-\frac{1}{2}\lambda_n, \alpha}},
\end{array} \]
   where
   $$ {\mathcal F}'_{\alpha,\beta}
     ={\mathcal F}^{1}_{\alpha,0}
      \otimes
      \bigotimes^n_{j=1}
      {\mbox Ker}_{{\mathcal F}^{2,j}_{0, l_j \vep_j}} Q^-_j , $$
   for $ \beta=l_1 \vep_1 + \cdots + l_n \vep_n $.
   Then we have the following proposition.

   \begin{proposition}(\cite{kn:JKM})
   Each ${\mathcal F}_i$ $(i=1,2,3,4)$ is
   an irreducible highest weight $U_q$-module
   isomorphic to $V(\mu_i)$,
   The highest weight vectors are given by
   $ | \mu_1 \rangle
    =|0, 0 \rangle $,
   $ | \mu_2 \rangle
    = b_1(-1) | \lambda_1, \lambda_1 \rangle $,
   $ | \mu_3 \rangle
    =| -\frac{1}{2}\lambda_n, 0 \rangle $,
   $ | \mu_4 \rangle
    =|-\frac{1}{2}\lambda_n-\vep_n, -\vep_n \rangle $.
   \end{proposition}


\section{Realization of the level one vertex operators}

We first recall the notion of Frenkel-Reshetikhin vertex operators
\cite{kn:FR}.
Let $V$ be a finite dimensional representation of the derived
quantum affine Lie
algebra $U_q'(\widehat{\mathfrak g})$
with the associated affinization space $V_z$
(cf. (\ref{S:level0}).
Let $V(\la)$ and $V(\mu)$ be two integrable highest weight
representations of $\U$.
The type I (resp. type II) {\it vertex operator} is the  $\U$-intertwining
operator $\tilde{\Phi}(z)$: $V(\la)\longrightarrow V(\mu)\hat{\otimes} V_z$
(resp.
$V(\la)\longrightarrow V_z\hat{\otimes} V(\mu)$), which equals to
 $\Phi(z)z^{\Delta_{\mu}-\Delta_{\la}}$ where $\Phi(z)$ is
the $U_q'(\widehat{\mathfrak g})$-intertwining operator.
Here the tensor product
$\hat{\otimes}$ is completed in certain sense, and we will write $\otimes$ for
$\hat{\otimes}$
in the sequel. $\Delta_{\mu}$ is the conformal weight for $\mu$.
For simplicity we will compute the
intertwining operators $\Phi(z)$ for the derived
subalgebra $U_q'(\widehat{\mathfrak g})$.

The existence of vertex operators
is given by the fundamental fact \cite{kn:FR}
(cf. \cite{kn:DJO}) (true for both types, though stated for type I):
\begin{eqnarray}\label{E:vertexpair}
Hom_{\U}(V(\lambda), V(\mu)\otimes V_z)
&\simeq &\{v\in V|\ wt(v)=\lambda-\mu \ \textstyle{mod} \ \delta
\ \ \mbox{and} \nonumber\\
&& e_i^{\langle \mu, h_i \rangle+1}v=0 \  \mbox {for} \ i=0, \cdots, n\},
\end{eqnarray}
where the isomorphism is defined as follows.
We say a pair of weights $(\la, \mu)$ is admissible if they satisfy
(\ref{E:vertexpair}).
For each admissible pair $(\la, \mu)$  there corresponds
uniquely a special vector $v_i$ in the crystal graph such that
$wt(v_i)=\la-\mu$. We then normalize the corresponding vertex operator
as follows.
\begin{equation} \label{E:4.3}
\Phi^{\mu V}_{\la}(z)|\la\rangle
=|\mu\rangle \otimes
v_{i}+\mbox{higher terms in $z$}
\end{equation}
Type II vertex operators assume similar normalization.

For the evaluation module $V_z$, we define the components of
the vertex operator
$\Phi_{\lambda}^{\mu V}:
V(\lambda)\longrightarrow V(\mu)\otimes V_z$ by
\begin{equation} \label{E:4.2}
\Phi_{\lambda}^{\mu V}(z)|u\rangle
=\sum_{j\in J}\Phi_{\lambda j}^{\mu V}
(z)|u\rangle \otimes v_j,
\end{equation}
for $|u \rangle \in V(\lambda)$, and $j$ runs through
the index set $J$
for the basis of $V$ as in the Figure 1.

We also consider the intertwining operators of modules of the following form:
$$
\Phi_{\lambda V}^{\mu}(z): V(\lambda)\otimes V_z
\longrightarrow V(\mu)\otimes {\bf C}[z, z^{-1}]
$$
by means of the vertex operators with respect to the dual space $V_z^*$:
\begin{equation}\label{E:4.1}
\Phi_{\lambda V}^{\mu}(z) (|v \rangle
\otimes v_i)=\Phi_{\lambda i}^{\mu V^*}(z)|v \rangle
\end{equation}
for $|v \rangle \in V(\lambda)$ and $i$ belongs to $J$.

For a node $j$ in the crystal graph we define the weight of
$j$ to be $\sum_i(f_i(j)-e_i(j))\Lambda_i$, where $f_i(j)$
is the number of $i$-arrows going out of $j$, and $e_i(j)$
is the number of $i$-arrows coming into $j$.

By the intertwining property it is easy to see the following
determination relations. For more explicit relations see \cite{kn:Ko1},
\cite{kn:JM2},
\cite{kn:JK}, \cite{kn:JKK}.

\begin{prop}\label{P:4.1}
The vertex operator $\Phi(z)$ of type I with respect to $V_z$ is
determined by its component $\Phi_{\overline 1}(z)$, where $\overline 1$
is the last vertex in the representation graph. In particular, for each
pair $i\stackrel{k}{\longrightarrow} j$ in the graph, we have
\begin{equation}
\Phi_i(z)=[r_j]_k^{-1}[\Phi_{j}(z), f_k]_{q_k^{-(wt(j), h_k)}},
\end{equation}
where $r_j$ counts the number of $k$-arrows coming into $j$.
With respect to $V_z^*$, the vertex operator $\Phi^*(z)$
of type I is determined by
$\Phi_1^*(z)$, and for each pair $i\stackrel{k}{\longrightarrow} j$
\begin{equation}
\Phi_{j}^*(z)=[r_j]_k^{-1}[f_k, \Phi^*_i(z)]_{q_k^{-(wt(i), h_k)}} ,
\end{equation}
where $r_j$ counts the number of $k$-arrows coming into $j$.
\end{prop}

\begin{prop}\label{P:4.2}
Let ${\Phi}(z)$ be a type II vertex operator with respect to $V_z$:
$V(\lambda)\rightarrow V_z\otimes V(\mu)$.
Then $ {\Phi}(z)$ is determined
by the component ${\Phi}_1(z)$. More precisely, for each with pair
$i\stackrel{k}{\longrightarrow} j$ we have:
\begin{equation}
\Phi_{j}(z)=[r_i']_k^{-1}[\Phi_i(z), e_k]_{q_k^{(wt(i), h_k)}}.
\end{equation}
where $r_i'$ counts the number of $k$-arrows going out of $i$.
With respect to $V_z^*$ the vertex operator $\Phi^*(z)$ of type II
is determined by its component
$\Phi_{\overline{1}}^*(z)$, where $\overline{1}$ is the last vertex
in the crystal graph. For each pair $i\stackrel{k}{\longrightarrow} j$
we have
\begin{equation}
\Phi_i^*(z)=q_k^2[r_i']_k^{-1}[e_k, \Phi_{j}^*(z)]_{q_k^{(wt(j), h_k)}} ,
\end{equation}
where $r_i'$ counts the number of $k$-arrows going out of $i$.
\end{prop}

\begin{prop} \label{P:4.3} Let $\U$ be the quantum affine algebra of type
$A_n^{(1)}, B_n^{(1)}$, or $D_n^{(1)}$.
(a) Let ${\Phi}(z): V(\lambda)\longrightarrow V(\mu)\otimes V_z$
be a vertex operator of type I,
where  $(\lambda, \mu)$ is an admissible pair of weights.
Then we have for each $j=1, \cdots, n$ and $k\in {\mathbb N}$
\begin{eqnarray*}
[{\Phi}_{\overline 1}(z), X_j^+(w)]&=&0, \quad
\mbox{for type A, we let $\overline 1=n+1$}\\
t_j{\Phi}_{\overline 1}(z)t_j^{-1} &=&
q^{\delta_{j1}}{\Phi}_{\overline 1}(z),\\
\ [a_j(k), {\Phi}_{\overline 1}(z)]&=&\delta_{j1}
\frac{[k]}{k} q^{(\tau+3/2)k}z^k{\Phi}_{\overline 1}(z),\\
\ [a_j(-k), {\Phi}_{\overline 1}(z)]&
=&\delta_{j1}\frac{[k]}{k} q^{-(\tau+1/2)k}z^{-k}
{\Phi}_{\overline 1}(z).
\end{eqnarray*}
(b) If ${\Phi}(z)$ is a vertex operator of type I
associated with $V_z^*$, then
\begin{eqnarray*}
[{\Phi}^*_1(z), X^+_j(w)]&=&0,
\quad\mbox{for type A, we let $\overline 1=n+1$}\\
t_j{\Phi}^*_1(z)t_j^{-1}&=&q^{\delta_{j1}}{\Phi}^*_1(z),\\
\ [a_j(k), {\Phi}^*_1(z)]&=&\delta_{j1}\frac{[k]}{k}
q^{\frac{3k}2}z^k{\Phi}^*_1(z),\\
\ [a_j(-k), {\Phi}^*_1(z)]&=&\delta_{j1}\frac{[k]}{k} q^{-\frac k2}
z^{-k}{\Phi}^*_1(z).
\end{eqnarray*}
(c) If ${\Phi}(z)$ is a vertex operator of type II
associated with $V_z$, then
\begin{eqnarray*}
[{\Phi}_1(z), X^-_j(w)]&=&0,\\
t_j{\Phi}_1(z)t_j^{-1}&=&q^{-\delta_{j1}}{\Phi}_1(z),\\
\ [a_j(k), {\Phi}_1(z)]&=&-\delta_{j1}\frac{[k]}{k}
q^{k/2}z^k{\Phi}_1(z),\\
\ [a_j(-k), {\Phi}_1(z)]&=&-\delta_{j1}\frac{[k]}{k}
q^{-3k/2}z^{-k}{\Phi}_1(z).
\end{eqnarray*}
(d) If ${\Phi}^*(z)$ is a vertex operator of type II associated
with $V_z^*$, then
\begin{eqnarray*}
[{\Phi}^*_{\overline 1}(z), X^-_j(w)]&=&0,\\
t_j{\Phi}^*_{\overline 1}(z)t_j^{-1}&=&q^{-\delta_{j1}}
{\Phi}^*_{\overline 1}(z),\\
\ [a_j(k), {\Phi}^*_{\overline 1}(z)]&=&-\delta_{j1}\frac{[k]}{k}
\left\{\begin{array}{ll}
q^{\frac{2n+3}2k}z^k{\Phi}^*_{\overline 1}(z), & \mbox{type $A_n^{(1)}$}\\
q^{-(\tau-1/2) k}z^k{\Phi}^*_{\overline 1}(z), & \mbox{type $B_n^{(1)}$ or
$D_n^{(1)}$}
\end{array}\right. \\
\ [a_j(-k), {\Phi}^*_{\overline 1}(z)]&=&-\delta_{j1}\frac{[k]}{k}
\left\{\begin{array}{ll}
q^{-\frac{2n+5}2 k}z^{-k}{\Phi}^*_{\overline 1}(z), & \mbox{type $A_n^{(1)}$}\\
 q^{(\tau-3/2)k}z^{-k}{\Phi}^*_{\overline 1}(z), &
\mbox{type $B_n^{(1)}$ or
$D_n^{(1)}$.}
\end{array}\right.
\end{eqnarray*}
\end{prop}

We introduce the auxiliary elements $a_{\overline 1}(k) \in
U_q(\hat{\bf h}^-)$ ($k\in {\mathbb Z}^{\times}$) such that
$$
[a_i(k), a_{\overline{1}}(l)]=\delta_{i1}\delta_{k, -l},
$$
and for type $A_n^{(1)}$ we need the elements $a_{n+1}(k)$
($k\in \mathbb Z^{\times}$) such that
$$
[a_i(k), a_{n+1}(l)]=\delta_{i,n}\delta_{k, -l}.
$$

>From Propositions (\ref{P:4.1}), (\ref{P:4.2}) and  (\ref{P:4.3})
we only need to
determine one component for each vertex operator.

\begin{theorem} The $\overline 1$-components of the type I vertex operator
${\Phi}(z)_{\lambda}^{\mu V}$ with respect to $V_z$
$: V(\lambda)\longrightarrow V(\mu)\otimes V_z$ can be realized explicitly as
follows:
\begin{eqnarray*}
{\Phi}_{\overline 1}(z)
&=exp(\sum\frac{[k]}kq^{(\tau +3/2)k}a_{\overline 1}(-k)z^k)
exp(\sum\frac{[k]}kq^{-(\tau+1/2)k}a_{\overline 1}(k)z^{-k})\\
&\hskip 1in \times e^{\la_1}(q^{\tau+1}z)^{\partial_{\la_1}+(\la_1|\la_1-
\overline{\mu})}(-1)^{\partial_{\la_n}+(\la_n|\la_1-\overline{\mu})}
b_{\la}^{\mu},
\end{eqnarray*}
for type $B_n^{(1)}$ and $D_n^{(1)}$, and where
$b_{\la}^{\mu}$ is a constant.
\begin{eqnarray*}
{\Phi}_{n+1}(z)
&=exp(\sum\frac{[k]}kq^{(n+5/2)k}a_{n+1}(-k)z^k)
exp(\sum\frac{[k]}kq^{-(n+3/2)k}a_{n+1}(k)z^{-k})\\
&\hskip 1in \times e^{\la_1}(q^{n+2}z)^{\partial_{\la_1}+(\la_1|\la_1-
\overline{\mu})}(-1)^{\partial_{\la_n}+(\la_n|\la_1-\overline{\mu})}
b_{\la}^{\mu},
\end{eqnarray*}
for type $A_n^{(1)}$ and where $b_{\la}^{\mu}$ is a constant.

The $1$-component of type I vertex operator $\Phi_{\la}^{\mu V^*}(z)$
associated with $V_z^*$ is given by
\begin{eqnarray*}
{\Phi}_1^*(z)&=&exp(\sum\frac{[k]}kq^{3k/2}a_{\overline 1}(-k)z^k)
exp(\sum\frac{[k]}kq^{-k/2}a_{\overline 1}(k)z^{-k})\\
\ && \hskip 1in\times e^{\la_1}(qz)^{\partial_{\la_1}+(\la_1|\la_1-
\overline{\mu})}(-1)^{\partial_{\la_n}+(\la_n|\la_1-\overline{\mu})}
b_{\la}^{\mu},
\end{eqnarray*}
where $b_{\la}^{\mu}$ is a constant.
\end{theorem}

\begin{theorem} The $1$-components of the type II vertex operator
${\Phi}_{\lambda}^{V\mu}(z)$ with respect to $V_z$
$: V(\lambda)\longrightarrow V_z\otimes V(\mu)$ can be realized explicitly as
follows:
\begin{eqnarray*}
{\Phi}_{1}(z)
&=&exp(-\sum\frac{[k]}kq^{\frac k2}a_{\overline 1}(-k)z^k)
exp(-\sum\frac{[k]}kq^{\frac{-3k}2}a_{\overline 1}(k)z^{-k})\\
\ &&\hskip 1in \times e^{-\la_1}(qz)^{-\partial_{\la_1}+(\la_1|\la_1+
\overline{\mu})}(-1)^{-\partial_{\la_n}+(\la_n|\la_1+\overline{\mu})}
b_{\la}^{\mu},
\end{eqnarray*}
where $b_{\la}^{\mu}$ is a constant.

The $\overline{1}$-component of type I vertex operator
$\Phi_{\overline{1}}^*(z)$
associated with $V_z^*$ is given by
\begin{eqnarray*}
{\Phi}_{\overline{1}}^*(z)
&=exp(-\sum\frac{[k]}kq^{-(\tau-1/2)k}a_{\overline 1}(-k)z^k)
exp(-\sum\frac{[k]}kq^{(\tau-3/2)k}a_{\overline 1}(k)z^{-k})\\
&\hskip 1in\times e^{-\la_1}(q^{-2n+2}z)^{-\partial_{\la_1}+(\la_1|\la_1+
\overline{\mu})}(-1)^{-\partial_{\la_n}+(\la_n|\la_1+\overline{\mu})}
b_{\la}^{\mu},
\end{eqnarray*}
for type $B_n^{(1)}$ and $D_n^{(1)}$, and where $b_{\la}^{\mu}$ is a constant.
\begin{eqnarray*}
{\Phi}_{n+1}^*(z)
&=exp(-\sum\frac{[k]}kq^{{\frac{2n+3}2}k}a_{n+1}(-k)z^k)
exp(-\sum\frac{[k]}kq^{-\frac{2n+5}2k}a_{n+1}(k)z^{-k})\\
&\hskip 1in\times e^{-\la_1}(q^{n+2}z)^{-\partial_{\la_1}+(\la_1|\la_1+
\overline{\mu})}(-1)^{-\partial_{\la_n}+(\la_n|\la_1+\overline{\mu})}
b_{\la}^{\mu},
\end{eqnarray*}
for type $A_n^{(1)}$ and where $b_{\la}^{\mu}$ is a constant.
\end{theorem}
\medskip

For type $C_n^{(1)}$ level $-1/2$ the vertex operators are determined by the
following result.

\begin{theorem} \cite{kn:JK}
   For the type one vertex operators associated to the admissible pair of
weights $(\lambda, \mu)$=
   $(\mu_1, \mu_2)$, $(\mu_2, \mu_1)$, $(\mu_3, \mu_4)$, and $(\mu_4, \mu_3)$
   respectively, one has
   \begin{eqnarray} \nonumber
      &&\widetilde{\Phi}^{\mu_j V}_{\mu_i}{}_{\overline 1}(z)
      =exp(\sum_{k=1}^{\infty}\frac{[k/2]}kq^{(n+1/4-\delta_{n1})k}
a_{\overline 1}(-k)z^k)\nonumber\\
&&\qquad\qquad
      exp(\sum_{k=1}^{\infty}\frac{[k/2]}kq^{-(n+3/4-\delta_{n1})k}
a_{\overline 1}(k)z^{-k}) \nonumber\\
    &&\qquad\qquad
  e^{\lambda_1}(q^{(n+1/2)}z)^{-\lambda_1(0)+1-(\lambda_1|\mu_i)}
\partial Z^+_1(q^{n+1/2}z)c_{ij},
\nonumber
 \end{eqnarray}
where $c_{ij}$ are constants
for the four cases $(\mu_i, \mu_j)$ with $c_{12}=1$.

 Moveover the type two vertex operators are given by
   \begin{eqnarray*}
      {\Phi}^{V\mu}_{\lambda}{}_1(z)
     &=&exp(-\sum_{k=1}^{\infty}\frac{[k/2]}kq^{-k/4}
a_{\overline 1}(-k)z^k)
      exp(-\sum_{k=1}^{\infty}\frac{[k/2]}kq^{-3k/4}
a_{\overline 1}(k)z^{-k})c_{ij}'\\
     &&\qquad\qquad
e^{-\lambda_1}(q^{-1/2}z)^{\lambda_1(0)+1+(\lambda_1|\mu_i)}
 Z^+_1(q^{-1/2}z)c_{ij}'
,
\end{eqnarray*}
where $c'_{ij}$ are constants
for the four cases $(\mu_i, \mu_j)$ with $c'_{12}=1$.
  \end{theorem}

\end{document}